\font\mathbb=msbm10
\newcount\counter\counter=0
\def\num{\global\advance \counter by 1 \the\counter}
\magnification=1100
\baselineskip=18pt
\parskip=6pt
\def\qed{\line{\hfill q.e.d.}}
\overfullrule=0pt

\font\eightpt=cmr8
\centerline{\bf COMMUTATORS IN THE ENDOMORPHISM RING OF A COMPLEX}
\centerline{\bf by}
\centerline{\bf Steven E. Landsburg}
\centerline{\bf University of Rochester}

\medskip

{\narrower\narrower\eightpt  According to [AM], the commutators in 
the endomorphism ring of a finite dimensional vector space are 
precisely the elements of trace zero.  We replace the finite 
dimensional vector space with a complex of finite dimensional 
vector spaces, and characterize commutators and other elements 
with commutator-like properties in terms of appropriately defined 
traces.\par\par}

\medskip
Let $V_\bullet$ be a complex of finite dimensional vector spaces 
over a field $k$.  We index cohomologically, so that the boundary 
map $d_i$ maps $V_i$ to $V_{i+1}$.  

{\bf Definition.}  $V_\bullet$ is {\sl quasi-bounded} if there 
exists at least one index $i$ for which $d_i:V_i\rightarrow 
V_{i+1}$ is the zero map.

For example, any complex that is bounded 
above or below is quasi-bounded.  

{\bf Definition.}  A {\it stretch} for $V_\bullet$ is  a finite 
interval $S=[M,N]\subset{\mathbb Z}$ that is maximal with respect 
to the property that for any $i,i+1\in S$, the map 
$d_i:V_i\rightarrow V_{i+1}$ is non-zero:

$$\matrix{
\bullet\bullet\bullet
&
\longrightarrow
&
V_{m-1}
&
\buildrel 0\over\longrightarrow&
\underbrace{\matrix{V_M
&
\longrightarrow&
\bullet\bullet\bullet&
\longrightarrow&
V_N}}_{\hbox{all maps non-zero in this range}}&
\buildrel 0\over\longrightarrow&
\longrightarrow&
\bullet\bullet\bullet\cr}$$

Note that a given complex might or might not have any stretches.

Now fix an endomorphism $\phi_\bullet:V_\bullet\rightarrow 
V_\bullet$.  Let $\phi_i^H:H^i(V_\bullet)\rightarrow 
H^i(V_\bullet)$ be the action of $\phi$ on cohomology.  

{\bf Definitions.} For an endomorphism $f$ of a vector space $V$, 
write $tr(f)$ for the trace of $f$.  Then for each index 
$i\in{\mathbb Z}$ and for each stretch $S\subset {\mathbb Z}$, 
set:

\itemitem{}$tr_i(\phi_\bullet)=tr(\phi_i)$
\itemitem{}$tr_i^H(\phi_\bullet)=tr(\phi^H_i)$
\itemitem{}$tr_S(\phi_\bullet)=\sum_{i\in S}(-1)^itr_i(\phi)$
\itemitem{}$tr_S^H(\phi_\bullet)=\sum_{i\in S}(-1)^itr_i^H(\phi)$

{\bf Definitions.} $\phi_\bullet$ is a {\sl commutator\/} if it is 
a commutator in the endomorphism ring of the complex $V_\bullet$.  
$\phi_\bullet$ is a {\sl pointwise commutator} if each 
$\phi_i:V_i\rightarrow V_i$ is a commutator in the endomorphism 
ring of the vector space $V_i$.

\medskip

We will prove the following theorems:

{\bf Theorem 1.}  The following are equivalent:
\itemitem{a)}$\phi_\bullet$ is a pointwise commutator.
\itemitem{b)}$tr_i(\phi_\bullet)=0$  for each index $i$.

{\bf Theorem 2.}  Suppose the field $k$ is infinite and 
$V_\bullet$ is quasi-bounded.  Then the following are equivalent:
\itemitem{a)}$\phi_\bullet$ is a commutator.
\itemitem{b)}$tr_i(\phi_\bullet)=tr_i^H(\phi_\bullet)=0$ for each 
index $i$.

{\bf Theorem 3.}  The following are equivalent:
\itemitem{a)}$\phi_\bullet$ is homotopic to a commutator.
\itemitem{b)}$tr_i^H(\phi_\bullet)=0$ for each index $i$.

{\bf Theorem 4.}  The following are equivalent:
\itemitem{a)}$\phi_\bullet$ is homotopic to a pointwise 
commutator.
\itemitem{b)}$tr_S(\phi_\bullet)=0$ for each stretch $S$.
\itemitem{c)}$tr_S^H(\phi_\bullet)=0$ for each stretch $S$.

\medskip

Theorem 1 is an instance of the main theorem in [AM], which states 
that any map of finite dimensional vector spaces is a commutator 
if and only if it has trace zero.  

To prove the remaining theorems, let $B_i=im(d_{i-1})\subset V_i$ 
and let $H_i=H^i(V_\bullet)$.  
Then up to isomorphism,
we can write $V_\bullet$ (depicted horizontally) and 
$\phi_\bullet$ (depicted vertically) in the form:

\def\Biggg#1{{\hbox{$\left#1\vbox to 24pt{}\right.$}}}

$$\matrix{
&&&&\pmatrix{0&0&1\cr0&0&0\cr0&0&0\cr}\cr
\bullet\bullet\bullet&
\longrightarrow &
&B_i\oplus H_i \oplus B_{i+1}&
\hbox to 1in{\rightarrowfill}
&
&B_{i+1}\oplus H_{i+1}\oplus B_{i+2}&
\longrightarrow&
\bullet\bullet\bullet\cr
&\vbox to .2in{}\cr
&&&
\Biggg\downarrow\llap{$\pmatrix{
\phi_i^B&g_i&h_i\cr 0&\phi_i^H& k_i\cr 0&0&\phi^B_{i+1}\cr}$\hbox 
to 
7pt{}}&&&
\Biggg\downarrow\rlap{$\pmatrix{
\phi_{i+1}^B&g_{i+1}&h_{i+1}\cr 0&\phi_{i+1}^H& 
k_{i+1}\cr 0&0&\phi^B_{i+2}\cr}$}\cr
\vbox to .2in{}\cr
\bullet\bullet\bullet&
\longrightarrow &
&B_i\oplus H_i \oplus B_{i+1}&
\hbox to 1in{\rightarrowfill}
&
&B_{i+1}\oplus H_{i+1}\oplus B_{i+2}&
\longrightarrow&
\bullet\bullet\bullet\cr
&&&&\pmatrix{0&0&1\cr0&0&0\cr0&0&0\cr}\cr
}\eqno(1)$$

{\bf Lemma A.}  For each $i\in \hbox{\mathbb Z}$, let $M_i$ and 
$N_i$ be traceless square matrices of sizes $m_i$ and $n_i$ over 
an infinite field $k$.  Then it is possible to choose square 
matrices $p_i,q_i,s_i,t_i$ that satisfy all of the following 
conditions for each $i$:
$$M_i=p_iq_i-q_ip_i\eqno(A1)$$
$$N_i=s_it_i-t_is_i\eqno(A2)$$
$$p_i\otimes I_{n_i}-I_{m_i}\otimes s_i^T\hbox{ is 
invertible.}\eqno(A3)$$ 
$$q_{i+1}\otimes I_{m_i}-I_{m_{i+1}}\otimes q_i^T\hbox{ is 
invertible.}\eqno(A4)$$
$$s_i\otimes I_{m_{i+1}}-I_{n_i}\otimes p_{i+1}^T\hbox{ is 
invertible.}\eqno(A5)$$
(Here $I_k$ is the identity matrix of size $k$, the empty $0$ by 
$0$ matrix is taken to be invertible by convention, and $M^T$ 
represents the transpose of $M$.)

{\bf Proof.}  By [AM], we can find $p_i,q_i,s_i,t_i$ satisfying 
(A1) and (A2).  

Now to satisfy (A3) and (A5), replace each $s_i$ by $s_i+\lambda_i 
I_{n_i}$, where $\lambda_i$ is not an eignenvalue of either 
$p_i\otimes I_{n_i}-I_{m_i}\otimes s_i^T$ or 
$ s_i\otimes I_{m_{i+1}}-I_{n_i}\otimes p_{i+1}^T$.
This is always possible by the infinitude of $k$, and it does not 
spoil (A1) or (A2).

To satisfy (A4), suppose it is already satisfied for $|i|<k$.  
Then add appropriate scalar multiples of the identity to $q_{-k-
1}$ and $q_{k+1}$ to satisfy (A4) for $|i|=k$, and continue by 
induction.

\qed

{\bf Proof of Theorem 2.}  Clearly a) implies b).  For the other 
direction, choose (by quasi-boundedness) an $N$ such that $B_N=0$ 
and therefore $\phi^B_N=0$.  Then b) gives
$$0=tr_{N-1}(\phi)=tr(\phi_{N-1}^B)+tr(\phi_{N-
1}^H)+tr(\phi_N^B)=tr(\phi_{N-1}^B)+0+0=tr(\phi_{N-1}^B)$$
$$0=tr_N(\phi)=tr(\phi_N^B)+tr(\phi_N^H)+tr(\phi_{N+1}^B)
=0+0+tr(\phi_{N+1}^B)=tr(\phi_{N+1}^B)$$
and, by backward and forward induction, 
$tr(\phi_i^B)=0$ for all $i$.

Therefore if we set $M_i=\phi_i^B$ and $N_i=\phi_i^H$, we can 
apply Lemma A to get matrices $p_i,q_i,s_i,t_i$ satisfying the 
conditions (A1) through (A5).

Now  to write $\phi_\bullet$ as a commutator, it suffices to write 
each $\phi_i$ as a commutator of the form
$$\phi_i=\pmatrix{\phi_i^B&g_i&h_i\cr 0&\phi_i^H&k_i\cr 
0&0&\phi^B_{i+1}\cr}=\left[\pmatrix{p_i&S_i&T_i\cr 0&s_i&U_i\cr 
0&0&p_{i+1}\cr},
\pmatrix{q_i&X_i&Y_i\cr 0&t_i&Z_i\cr 0&0&q_{i+1}\cr}\right]$$

I claim that such an expression exists with $S_i=U_i=Y_i=0$.  To 
establish the claim it suffices to solve the equations
$$p_iX_i-X_is_i=g_i\qquad T_iq_{i+1}-q_iT_i=h_i\qquad s_iZ_i-
Z_ip_{i+1}=k_i\eqno(3)$$
for the unknown matrices $X_i, T_i, Z_i$.
This  solvability follows from (A3), (A4) and (A5) of Lemma A.

\qed 

{\bf Proof of Theorem 3.}  a) implies b) immediately.  For the 
other direction,  
represent $\phi$ as in (1).  Define $S:V_i\rightarrow V_{i-1}$ by 
the matrix
$$S_i=\pmatrix{0&0&0\cr
k_{i-1}&0&0\cr
\phi_i^B&g_i&h_i\cr}$$
  
Then $S_\bullet$ is a homotopy from $\phi$ to the map represented 
by
$$\pmatrix{0&0&0\cr 0&\phi_i^H&0 \cr 0&0&0\cr}\eqno(4)$$
so we can assume $\phi$ has the form (4).  Now write the traceless 
map of vector spaces $\phi_i^H$ as a commutator $[s_i,t_i]$ and 
note that 
$$\phi_\bullet=\left[\pmatrix{0&0&0\cr 0&s_\bullet&0\cr 0&0&0\cr},
\pmatrix{0&0&0\cr 0&t_\bullet&0\cr 0&0&0\cr}\right]$$

\qed

{\bf Lemma B.}  Let $\{T_i|i\in{\mathbb Z}\}$ be a family of 
scalars 
with the property that for any stretch $S$ of $V_\bullet$, we have 
$\sum_{i\in S}(-1)^iT_i=0$.  Then there exists a null-homotopic 
map $\tau:V_\bullet\rightarrow V_\bullet$ such that $tr(t_i)=T_i$.  

{\bf Proof.}  First I claim that we can find scalars $S_i\in k$ 
satisfying the two conditions
\itemitem{a)}$T_i=S_i+S_{i+1}$ for each $i$.
\itemitem{b)}Whenever $B_i=0$, $S_i=0$.  (Recall that $B_i$ is the 
image of the map $\phi_{i-1}:V_{i-1}\rightarrow V_i$.)

To establish the claim, choose any index $n$ with the property 
that $B_i=0$, if there is any such index, or choose $n$ completely 
arbitrarily otherwise.  Define $S_n=0$ and use a) with forward and 
backward induction to define all other $S_i$.  The hypothesis of 
the lemma then guarantees that b) also holds, which proves the 
claim. 

Now let $\tau_i^B:B_i\rightarrow B_i$ be any map with trace $S_i$, 
and define the map $\tau$ as follows:

$$\matrix{
&&&&\pmatrix{0&0&1\cr0&0&0\cr0&0&0\cr}\cr
\bullet\bullet\bullet&
\longrightarrow &
&B_i\oplus H_i \oplus B_{i+1}&
\hbox to 1in{\rightarrowfill}
&
&B_{i+1}\oplus H_{i+1}\oplus B_{i+2}&
\longrightarrow&
\bullet\bullet\bullet\cr
&\vbox to .2in{}\cr
&&&
\Biggg\downarrow\llap{$\pmatrix{
\tau_i^B&0&0\cr 0&0& 0\cr 0&0&\tau_{i+1}^B\cr}$\hbox to
7pt{}}&&&
\Biggg\downarrow\rlap{$\pmatrix{
\tau_{i+1}^B&0&0\cr 0&0& 0\cr 0&0&\tau_{i+2}^B\cr}$}\cr
\vbox to .2in{}\cr
\bullet\bullet\bullet&
\longrightarrow &
&B_i\oplus H_i \oplus B_{i+1}&
\hbox to 1in{\rightarrowfill}
&
&B_{i+1}\oplus H_{i+1}\oplus B_{i+2}&
\longrightarrow&
\bullet\bullet\bullet\cr
&&&&\pmatrix{0&0&1\cr0&0&0\cr0&0&0\cr}\cr
}$$

\qed

{\bf Proof of Theorem 4.}  Clearly b) and c) are equivalent, and 
clearly a) implies both.  For the other direction, assume b) and 
put $T_i=tr_i(\phi)$.
Then after subtracting the null-homotopic map 
$\tau$ constructed in Lemma B, we can assume $tr_i(\phi)=0$ for 
all 
$i$, and the result now follows from Theorem 1.

\qed

\vfill\eject

\centerline{\bf Counterexamples}

Theorem 2, the proof of which depends on Lemma A, assumes quasi-
boundedness and an infinite base field.
Example 1 below will show that Theorem 2 can fail in the absence 
of quasi-boundedness.  Example 2 below will show that Lemma A (but 
not necessarily Theorem 2) can fail in the absence of an infinite 
field.

{\bf Example 1.}  
For all $i$, let $V_i=k\oplus k$, let $d_i:V_i\rightarrow V_{i+1}$ 
be the map $\pmatrix{0&1\cr 0&0}$, and let $\phi_i:V_i\rightarrow 
V_i$ be the map $(-1)^i\pmatrix{1&0\cr 0&-1}$.  

Then $\phi_\bullet$ is a map of complexes satisfying condition b) 
of Theorem 2.  But $\phi_\bullet$ cannot be a commutator, because 
if it were, the induced map on boundaries would also be a 
commutator and therefore pointwise traceless, whereas in fact that 
induced map has trace $\pm 1$.  Thus Theorem 2 can fail when 
$V_\bullet$ is not quasi-bounded.

{\bf Example 2.}  Let $k=\hbox{\mathbb F}_2$.  For all $i$, let
$$M_i=\pmatrix{0&0\cr 1&0\cr}$$
$$N_i=\pmatrix{0&1\cr 0&0\cr}$$

If Lemma A holds over ${\mathbb F}_2$, 
then we must be able to write
$$M_1=p_1q_1-q_1p_1\qquad N_1=s_1t_1-t_1s_1
\qquad M_2=p_2q_2-q_2p_2 \eqno(5)$$
so that the following matrices are invertible:

$$p_1\otimes I_2-I_2\otimes s_1^T\eqno(6)$$
$$ q_2\otimes I_2-I_2\otimes q_1^T\eqno(7)$$
$$ p_2\otimes I_2-I_2\otimes s_1^T\eqno(8)$$

For each matrix $M$, put $C(M)=\{P|\exists Q \hbox{ such that } 
PQ-QP=M\}$.  Then (5) requires that:

\def\mattwo{\pmatrix{0&0\cr0&1\cr}}
\def\matthree{\pmatrix{0&0\cr1&0\cr}}
\def\matfour{\pmatrix{0&0\cr1&1\cr}}
\def\matfive{\pmatrix{0&1\cr0&0\cr}}
\def\matsix{\pmatrix{0&1\cr0&1\cr}}

\def\matnine{\pmatrix{1&0\cr0&0\cr}}
\def\mateleven{\pmatrix{1&0\cr1&0\cr}}
\def\mattwelve{\pmatrix{1&0\cr1&1\cr}}
\def\matthirteen{\pmatrix{1&1\cr0&0\cr}}
\def\matfourteen{\pmatrix{1&1\cr0&1\cr}}

$$p_1,p_2\in
C(\phi_1^B)=C(\phi_2^B)=\left\{\mattwo,\matthree,\matfour,
\matnine,\mateleven,\mattwelve\right\}\eqno(9)$$
$$s_1\in
C(\phi_1^H)=\left\{\mattwo,\matfive,\matsix,\matnine,
\matthirteen,\matfourteen\right\}\eqno(10)$$

Checking each possible pair, the invertibility of (6) and (8) 
requires that 
$$(p_1,s_1),(p_2,s_1)\in\left\{\left(\matthree,\matfourteen\right)
,
\left(\mattwelve,\matfive\right)\right\}$$

In particular, for each $i=1,2$, either $p_i=\matthree$ or 
$p_i=\mattwelve$, either of which, in conjunction with (5),  
implies
that $$q_i\in\left\{\mattwo,\matfour,\matnine,\mateleven\right\}$$

Taking the sixteen possible ordered pairs for $(q_1,q_2)$, one 
checks 
that none of them makes (7) invertible. Thus Lemma A, which is used in
the proof of Theorem 2, can fail in the absence of the assumption that
$k$ is infinite.

Note that Lemma A is sufficient but not necessary for the solvability 
of the system (3), which in turn is sufficient but not necessary for the
proof of Theorem 2.

\medskip

\centerline{\bf References}

\noindent [AM] A.~Albert and B.~Muckenhoupt, {\it On Matrices of
Trace Zero}, Michican Math Journal 4 (1957).

\bye